\documentclass[reqno,12pt]{amsart}

\usepackage{epsf}
\usepackage{graphics}
\usepackage{graphicx}
\usepackage{amssymb}
\usepackage{amsmath}
\usepackage{hyperref}

\date{}

\theoremstyle{plain}
\newtheorem{theorem}{Theorem}

\newtheorem{rem}{Remark}

\theoremstyle{definition}

\theoremstyle{remark}

\def\R{{\mathbb R}}

\title{Systoles on Punctured Spheres}

\author{Sebastian Baader and Jasmin J\"org}

\begin{document}

\begin{abstract} We determine the maximal number of systoles among all spheres with $n$ punctures endowed with a complete Riemannian metric of finite area.
\end{abstract}


\maketitle

\section{Introduction}

Every standard orientable surface $\Sigma_{g,n}$ of genus~$g$ with $n$ punctures and strictly negative Euler characteristic, i.e.\ $2-2g-n<0$, admits a complete hyperbolic metric. It is an open problem to determine how many systoles such a hyperbolic surface can have. Here a systole is an essential closed geodesic of minimal length, where essential means it is neither homotopic to a point nor to a puncture. In the case of surfaces $\Sigma_{g,n}$ with $n \leq g$, the answer is known to be a function of $g$ roughly bounded between $g^{\frac{4}{3}-\epsilon}$ and $g^2/\log(g)$, for all $\epsilon>0$. The lower and upper bound are due to Schmutz Schaller and Fanoni-Parlier, respectively~\cite{Sch2,FP}. In the case of punctured spheres $\Sigma_{0,n}$, the answer is linear in $n$, more precisely between $n$ and $\frac{7}{2}n$, thanks to the very last statement in~\cite{FP}.

In this note, we derive the maximal number of systoles, known as kissing number, for punctured spheres $\Sigma_{0,n}$ with Riemannian metrics of arbitrary curvature. We denote by $\mathcal{RM}_{0,n}$ the space of complete, finite area Riemannian metrics on $\Sigma_{0,n}$. In the case of variable curvature, a given closed curve can be homotopic to several geodesics of shortest length. For that reason, we consider the set of systoles $\text{Sys}(S)$ of a surface $S$ up to homotopy, and define its kissing number as
$$\text{Kiss}(S)=|\text{Sys}(S)/\sim|,$$
where the equivalence relation $\sim$ denotes free homotopy between closed curves.

\begin{theorem}
\label{kiss}
For all $n \geq 5$,
$$\max\{\text{Kiss}(S) \mid S \in \mathcal{RM}_{0,n}\}=4n-11.$$
\end{theorem}

In simple words, the maximal kissing number of a punctured sphere is about four times the number of punctures. The condition $n \geq 5$ is necessary since the four-punctured sphere admits at most three systoles \cite{Sch1}.

As we will see in the next section, the upper bound follows easily from the work of Fanoni and Parlier~\cite{FP}. The lower bound, presented in the third section, is based on a construction, also inspired by~\cite{FP}, and a remark by Adams~\cite{A}.

\section{Upper bound on the kissing number}

We begin by recalling a few results about the kissing number of hyperbolic surfaces. For closed surfaces, Parlier derives an upper bound of order $g^2 / \log(g)$ in terms of the genus $g$ of the surface in~\cite{Pa}.
This is generalised to surfaces with $n$ punctures by Fanoni and Parlier in~\cite{FP}, where they provide an upper bound of order $(g+n)g/\log(g+1)$. Both articles make use of bounds on the maximal systole length. This quantity has been studied mostly for hyperbolic surfaces, e.g.\ Buser and Sarnak provided the first family of closed surfaces with systole length at least $4/3 \log(g)$ \cite{BuS}. This establishes a logarithmic order of growth, as a straightforward area argument yields a $2 \log(g)$ upper bound.

It is well known that on closed surfaces with any Riemannian metric, systoles intersect at most once. An upper bound for the kissing number in this more general setting is thus given by the maximal number of curves intersecting pairwise at most once. Greene gives an upper bound of order $g^2\log(g)$, in terms of the genus $g$ of the surface, for the maximal size of such curve systems \cite{G}. The intersection properties of systoles turn out to be a useful tool when considering punctured spheres; Fanoni and Parlier make extensive use of them in~\cite{FP}.

For the upper bound of Theorem~\ref{kiss}, we fix $S \in \mathcal{RM}_{0,n}$ ($n \geq 5$) and follow the proof given for Theorem~4.13 in~\cite{FP}, adapting slightly to the arbitrary curvature setting where necessary. The main ingredients are Propositions~3.2 and 3.3 in~\cite{FP}, stated below as \eqref{in} and \eqref{cu}, respectively. 
We briefly summarise and point out that the proofs given in that article are based on the fact that smoothing intersection points and rounding corners reduce length, so they work perfectly fine for Riemannian surfaces with arbitrary curvature.

\begin{enumerate}
    \item Any two systoles on $S$ intersect at most twice. \label{in}
    \item If two systoles on $S$ intersect twice, then one of them bounds two cusps. \label{cu}
\end{enumerate}

The proof of Theorem~4.13 in~\cite{FP} makes use of the fact that all closed curves on the punctured sphere $S$ are separating. As a consequence, systoles come in two types: ones that bound a pair of cusps on one side and ones that bound more than two cusps on either side. Following~\cite{FP}, we define
\begin{align*}
    A(S) = \{\alpha \in \text{Sys}(S) \, \mid \alpha \text{ bounds two cusps}\},
\end{align*}
and estimate $|A(S)|$ and $|\text{Sys}(S) \setminus A(S)|$ separately.
We start with the latter. Based solely on \eqref{in} and \eqref{cu}, the authors argue that all curves of $\text{Sys}(S) \setminus A(S)$ are pairwise disjoint. This implies
$$|\text{Sys}(S) \setminus A(S)| \leq n-5.$$
As the argument makes no use of curvature, it carries over to yield
$$|\text{Sys}(S) \setminus A(S)/\sim| \leq n-5$$
in our setting.

In contrast, their upper bound for the first set, $|A(S)| \leq \frac{5}{2}n$, does not carry over to the setting of variable curvature. However, they provide a simple topological proof of the weaker bound $|A(S)| \leq 3n-6$ in Remark~4.3. We give a brief sketch of the argument. Every systole in $A(S)$ bounds a disc with two punctures. In that disc, there is a unique arc connecting the two cusps, up to homotopy. The set of arcs obtained in that way can be chosen to have pairwise disjoint interiors (for example, by choosing them to be geodesic arcs with respect to an auxiliary hyperbolic metric on $S$). As a consequence, these arcs can be completed to a triangulation of $S$ with $n$ vertices, $e \geq |A(S)/\sim|$ edges, and $\frac{2}{3} e$ triangles. Using the Euler characteristic of the sphere, we obtain $2=n-\frac{1}{3} e$, hence
\begin{align*}
    |A(S)/\sim| \leq 3(n-2).
\end{align*}
Combining the bounds on $|A(S)/\sim|$ and $|\text{Sys}(S) \setminus A(S)/\sim|$ we get the desired upper bound for Theorem~\ref{kiss}:
$$ \text{Kiss}(S) \leq 4n-11.$$

\section{Lower bound on the kissing number}

In this section, we construct a Riemannian metric with at least ${4n-11}$ systoles on the sphere with $n$ punctures. Our construction is inspired by a construction of a hyperbolic plane with infinitely many cusps given in~\cite{A}, which we briefly recall. The building block is a particular neighbourhood of a hyperbolic cusp. We glue $12$ copies of the semi-ideal hyperbolic triangle with angles $0,\frac{\pi}{2},\frac{\pi}{3}$ in a suitable way, by successive mirroring along infinite edges, to obtain a neighbourhood of a hyperbolic cusp with area $2 \pi$, whose boundary is a regular hexagon with angles $\frac{2\pi}{3}$. The latter can be glued together in the pattern of the regular hexagonal lattice in $\R^2$. The resulting surface is a plane with countably many cusps and constant curvature $-1$. This construction also appears in~\cite{FP}. One of its features is that every pair of adjacent cusps is surrounded by a systole, so every cusp is surrounded by six systoles with pairwise intersection number two, see Section~3 in~\cite{FP} for more details.

In the following, we construct a discretised version of the above plane with cusps, which will serve to construct the desired sphere with $n+4$ punctures.
Consider the lattice $\Gamma \subset \R^2$ spanned by the two vectors $v_1=2e_1$, $v_2=e_1+\sqrt{3}e_2$. Observe that this lattice corresponds to a tiling of the plane with regular triangles of side length two. Remove a disc of radius $r<1$ around each lattice point, and glue in a Euclidean half-cylinder with base circle of length $2 \pi r$ to each arising boundary circle. The resulting surface $\widetilde{X}$ admits an isometric embedding into $\R^3$ as a horizontal plane with vertical `chimneys', one for each point of the lattice $\Gamma$. Of course, the surface $\widetilde{X}$ fails to be smooth and has infinite area. Note that the geodesics encircling the chimneys are not systoles, since they are homotopic to a puncture. 
The actual systoles of $\widetilde{X}$ are precisely the geodesics that wind around two neighbouring chimneys, meaning around chimneys whose base points in $\Gamma$ have distance two. Their length is $4+2 \pi r$; every other closed geodesic is strictly longer since it surrounds at least three chimneys or two chimneys further apart.

\begin{figure}[htb]
\includegraphics[scale=0.7]{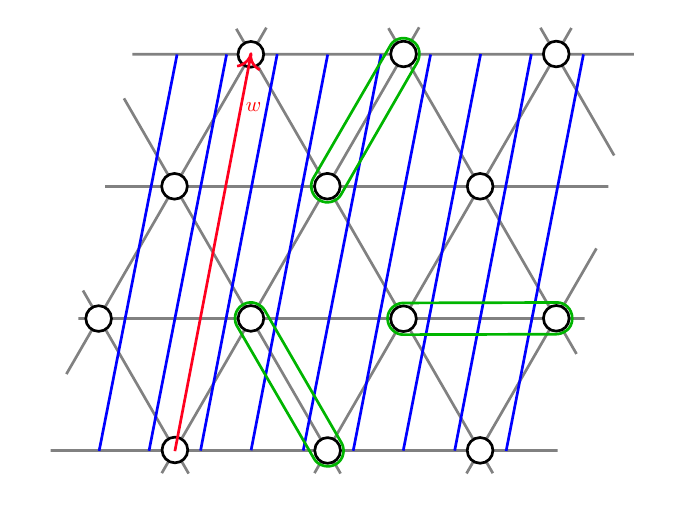}
\caption{Two types of systoles}
\label{meridian}
\end{figure}

On the surface $\widetilde{X}$ just constructed, there is a translation $T_w \colon \widetilde{X} \to \widetilde{X}$ by the vector $w=3v_2-v_1=e_1+3\sqrt{3}e_2$. We define $X$ to be the quotient surface of $\widetilde{X}$ by the isometry $T_w$: $X=\widetilde{X} / T_w$. We may view $X$ as a cylinder with three countable rows of chimneys. Forgetting the chimneys for a while, the length of the meridional geodesics of the cylinder is
$$m=\sqrt{(3\sqrt{3})^2+1^2}=\sqrt{28}.$$
Now, back to the chimneys, we choose their radius~$r$ so that the systoles surrounding neighbouring pairs match the length of the meridional geodesics: $4+2 \pi r=\sqrt{28}$, i.e. 
$$r=\frac{\sqrt{28}-4}{2 \pi} < \frac{1}{4}.$$
A careful inspection of Figure~\ref{meridian} reveals that the diameter of the chimneys is small enough so that some of the meridional geodesics of $X$ still run straight between them. More precisely, this holds because of the following: In the right-angled triangle with side lengths $\sqrt{3},\frac{1}{3},\frac{\sqrt{28}}{3}$ with hypotenuse $\frac{\sqrt{28}}{3}$, the height $h=\frac{\sqrt{3}}{\sqrt{28}} > \frac{1}{4}$ is strictly larger than $r$.

\begin{figure}[htb]
\includegraphics[scale=0.7]{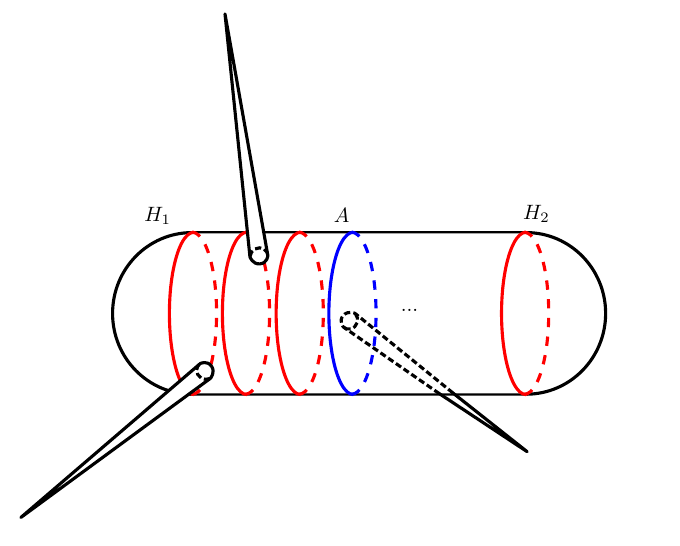}
\caption{A sphere $S$ with $n$ punctures}
\label{sphere}
\end{figure}

Now comes the final step of the construction: choose two meridional geodesics of $X$, which together bound a cylinder $A \subset X$ that contains $n$ cusps. Glue round hemispheres $H_1,H_2$ of circumference $\sqrt{28}$ to both boundary geodesics of $A$.
The resulting surface
$$S=A \cup H_1 \cup H_2$$
is a sphere with $n$ punctures, see Figure~\ref{sphere} for a schematic picture. A careful edge count reveals that the cylinder $A \subset S$ contains at least $3n-6$ systoles surrounding pairs of adjacent chimneys (the precise number is $0,1$ for $n=1,2$, respectively, and $3n-6$ for $n \geq 3$; in order to understand the leading term $3n$, recall that every point of the lattice $\Gamma$ has six nearest neighbours). As for the meridional geodesics, it appears that there are $n+1$ in number, up to homotopy. However, not all of them count as systoles, the ones not counting are marked red in Figure~\ref{sphere}. Firstly, the boundary geodesics become contractible after gluing in the hemispheres. Secondly, the meridians between the first and last pair of cusps surround a single cusps, so they do not count either. Thirdly, there are two further meridians that bound a disc with two cusps (the first two and the last two cusps, respectively). After an obvious homotopy, these appear to be in the set of $3n-6$ systoles surrounding pairs of adjacent chimneys, so we have already counted them. We are left with $n-5$ meridional systoles, up to homotopy. Here we use the assumption $n \geq 5$. By construction, all systoles of $S$ are homotopic to systoles contained in the annular part $A$. In order to see that the part $A$ contains no shorter closed geodesics, we remark that the ones that lift to closed geodesics in $\widetilde{X}$ must have length at least $\sqrt{28}$; the ones that do not lift to closed geodesics in $\widetilde{X}$ wind non-trivially around the cylinder $A$, thus also have length at least~$\sqrt{28}$.

In conclusion, the surface $S$ is indeed a sphere with $n$ punctures and $(3n-6)+(n-5)=4n-11$ systoles. However, $S$ fails to be smooth and has infinite area. This is easy to fix, in one stroke. We adapt the construction of $\widetilde{X}$, as follows: instead of replacing the discs around lattice points of $\Gamma \subset \R^2$ by vertical half-cylinders, we glue in vertical cusps with negative curvature that rise from the plane with horizontal tangent planes, much like a standing tuba. We may choose the curvature to be constant in a disc neighbourhood of the cusp, thus ensuring finite area. In order to keep the systole length to be $\sqrt{28}$, we must slightly enlarge the radius $r$, for example to $r=\frac{1}{4}$, and curve up the cusps fast enough so that the systoles surrounding two neighbouring cusps keep their length $\sqrt{28}$. Here we simply appeal to the intermediate value theorem. With this modification, the surface $\widetilde{X}$ becomes smooth. Of course, we can also choose a smooth transition between the annular part and hemispheres $H_{1,2}$. The resulting surface $S$ is then smooth, complete and with finite area.

\section{Final Remarks}

The surface $S$ constructed in the previous section has parts with negative and positive curvature. According to Fanoni and Parlier's last theorem in~\cite{FP}, there is no way of improving $S$ to an example with constant curvature $-1$ since, in that case, the kissing number is bounded above by $\frac{7}{2}n-5$. However, at the cost of a small constant drop in the kissing number, we can modify the surface $S$ to have non-positive curvature. Indeed, instead of capping off the annular part $A$ with two hemispheres in the above construction, we may cap it off with two hyperbolic pairs of pants $P_1,P_2$ with two cusps and one boundary geodesic of length $\sqrt{28}$.

A comment about the translation vector $w=e_1+3\sqrt{3}e_2$ is in order. The seemingly arbitrary choice of $w$ in our construction is in fact quite forced. Indeed, the analogous construction does not work with a number of rows of chimneys other than three. If we choose two rows of chimneys on the cylinder $X$, then the meridional geodesic becomes too short; if we choose at least four rows, then it becomes too long, independent of the choice of the radius $r$ of the chimneys.

We conclude by evoking a picture of the surface $S$. Despite the abstract nature of its construction, the surface $S$ admits a fairly adequate description as the surface of the stem and branches of a common spruce with constant circumference and three rows of equidistant branches.

\bigskip
\noindent
Mathematisches Institut, Universit\"at Bern, Sidlerstrasse 5, 3012 Bern, Switzerland

\smallskip
\noindent
\texttt{sebastian.baader@unibe.ch}

\smallskip
\noindent
\texttt{jasmin.joerg@unibe.ch}

\end{document}